\documentclass[10pt,leqno]{article}

\usepackage{amsmath,amssymb,amsthm,mathrsfs,dsfont}

\usepackage[margin=3cm]{geometry} 

\usepackage{titlesec,hyperref}
\usepackage{extarrows}

\usepackage{color}

\usepackage{fancyhdr}
\titleformat{\subsection}{\it}{\thesubsection.\enspace}{1pt}{}

\newtheorem{theo}{Theorem}[section]
\newtheorem{lemm}[theo]{Lemma}
\newtheorem{defi}[theo]{Definition}

\newtheorem{rema}[theo]{Remark}
\numberwithin{equation}{section}

\allowdisplaybreaks 

\begin{document}
\title{Global weak solutions for a three-component Camassa-Holm system with N-peakon solutions
\hspace{-4mm}
}
\author{Wei $\mbox{Luo}^1$\footnote{E-mail:  luowei23@mail2.sysu.edu.cn} \quad and\quad
 Zhaoyang $\mbox{Yin}^{1,2}$\footnote{E-mail: mcsyzy@mail.sysu.edu.cn}\\
 $^1\mbox{Department}$ of Mathematics,
Sun Yat-sen University,\\ Guangzhou, 510275, China\\
$^2\mbox{Faculty}$ of Information Technology,\\ Macau University of Science and Technology, Macau, China}

\date{}
\maketitle
\hrule

\begin{abstract}
In this paper we mainly investigate the Cauchy problem of a three-component Camassa-Holm system. By using the method of approximation of smooth solutions, a regularization technique and the special structure of the system, we prove the existence of global weak solutions to the system.\\

\vspace*{5pt}
\noindent {\it 2010 Mathematics Subject Classification}: 35Q53 (35B30 35B44 35C07 35G25)

\vspace*{5pt}
\noindent{\it Keywords}: A three-component Camassa-Holm system; the method of approximation; the regularization technique, global weak solutions.
\end{abstract}

\vspace*{10pt}

\tableofcontents
\section{Introduction}
  In this paper we consider the Cauchy problem for the following three-component Camassa-Holm equations with N-peakon solutions:
  \begin{align}
\left\{
\begin{array}{ll}
u_{t}=-va_x+u_xb+\frac{3}{2}ub_x-\frac{3}{2}u(a_xc_x-ac),\\[1ex]
v_t=2vb_x+v_xb,\\[1ex]
w_t=-vc_x+w_xb+\frac{3}{2}wb_x+\frac{3}{2}w(a_xc_x-ac),\\[1ex]
u=a-a_{xx},\\[1ex]
v=\frac{1}{2}(b_{xx}-4b+a_{xx}c_x-c_{xx}a_x+3a_xc-3ac_x), \\[1ex]
w=c-c_{xx},\\[1ex]
u|_{t=0}=u_{0},~~ v|_{t=0}=v_0,~~ w|_{t=0}=w_0.\\[1ex]
\end{array}
\right.
\end{align}
This system was proposed by Geng and Xue in \cite{Geng.Xue}. It is based on the following spectral problem
\begin{align}
\phi_x=U\phi,~~\phi=
\begin{pmatrix}
\phi_1 \\
\phi_2 \\
\phi_3
\end{pmatrix},~~~
U=
\begin{pmatrix}
0 & 1 & 0 \\
1+\lambda v & 0 & u \\
\lambda w   & 0 & 0
\end{pmatrix},
\end{align}
where $u,~v,~w$ are three potentials and $\lambda$ is a constant spectral parameter. It was shown in \cite{Geng.Xue} that the N-peakon solitons of the system (1.1) have the form
\begin{align}
a(t,x)=\sum^N_{i=0} a_i(t)e^{-|x-x_i(t)|},\\
\nonumber b(t,x)=\sum^N_{i=0} b_i(t)e^{-2|x-x_i(t)|},\\
\nonumber c(t,x)=\sum^N_{i=0} c_i(t)e^{-|x-x_i(t)|},
\end{align}
where $a_i,~b_i,~c_i$ and $x_i$ evolve according to a dynamical system. Moreover, the author derived infinitely many conservation laws of the system (1.1). By setting  $a=c=0$ , the system (1.1) reduces to
\begin{align}
v_t=2v_xb+vb_x, ~~v=\frac{1}{2}(b_{xx}-4b).
\end{align}
Taking advantage of an appropriate scaling $\widetilde{v}(t,x)=v(\frac{t}{2},\frac{x}{2}),~~\widetilde{b}(t,x)=-b(t,\frac{x}{2})$, one can deduce that
\begin{align}
\widetilde{v}_t+2\widetilde{v}_x\widetilde{b}+\widetilde{v}\widetilde{b}_x=0, ~~\widetilde{v}=\widetilde{b}-\widetilde{b}_{xx},
\end{align}
which is nothing but the famous Camassa-Holm (CH) equation \cite{Camassa, Constantin.Lannes}. The Camassa-Holm equation was derived as a model for shallow water waves \cite{Camassa, Constantin.Lannes}. It has been investigated extensively because of its great physical significance in the past two decades. The CH equation has a bi-Hamiltonian structure \cite{Constantin-E,Fokas} and is completely integrable \cite{Camassa,Constantin-P}. The solitary wave solutions of the CH equation were considered in \cite{Camassa,Camassa.Hyman}, where the authors showed that the CH equation possesses peakon solutions of the form $Ce^{-|x-Ct|}$. It is worth mentioning that the peakons are solitons and their shape is alike that of the travelling water waves of greatest height, arising as solutions to the free-boundary problem for incompressible Euler equations over a flat bed (these being the governing equations for water waves),
cf. the discussions in \cite{Constantin2,Constantin.Escher4,Constantin.Escher5,Toland}. Constantin and Strauss verified that the peakon solutions of the CH equation are orbitally stable in \cite{Constantin.Strauss}. \\
  $~~~~~~$The local well-posedness for the CH equation was studied in \cite{Constantin.Escher,Constantin.Escher2,Danchin,Guillermo}. Concretely, for initial profiles $\widetilde{b}_0\in H^s(\mathbb{R})$ with $s>\frac{3}{2}$, it was shown in \cite{Constantin.Escher,Constantin.Escher2,Guillermo} that the CH equation has a unique solution in $C([0,T);H^s(\mathbb{R}))$. Moveover, the local well-posedness for the CH equation in Besov spaces $C([0,T);B^s_{p,r}(\mathbb{R}))$ with $s>\max(\frac{3}{2},1+\frac{1}{p})$ was proved in \cite{Danchin}. The global existence of strong solutions were established in \cite{Constantin,Constantin.Escher,Constantin.Escher2} under some sign conditions and it was shown in \cite{Constantin,Constantin.Escher,Constantin.Escher2,Constantin.Escher3} that the solutions will blow up in finite time when the slope of initial data was bounded by a negative quantity. The global weak solutions for the CH equation were studied in \cite{Constantin.Molinet} and \cite{Xin.Z.P}. The global conservative and dissipative solutions of CH equation were presented in \cite{Bressan.Constantin} and  \cite{Bressan.Constantin2}, respectively.\\
$~~~~~~$A natural idea is to extend such study to the multi-component generalized
systems. One of the most popular generalized systems is the following integrable two-component
Camassa-Holm shallow water system (2CH) \cite{Constantin.Ivanov}:
\begin{align}
\left\{
\begin{array}{ll}
m_t+um_x+2u_xm+\sigma\rho\rho_x=0, \\[1ex]
\rho_t+(u\rho)_x=0,
\end{array}
\right.
\end{align}
where $m=u-u_{xx}$ and $\sigma=\pm1$. Local well-posedness for (2CH) with the initial
data in Sobolev spaces and in Besov spaces was established in \cite{Constantin.Ivanov}, \cite{Escher.Yin}, and  \cite{GuiGuilong}, respectively. The blow-up phenomena and global existence of strong solutions to (2CH) in
Sobolev spaces were obtained in \cite{Escher.Yin}, \cite{Guan.Yin} and \cite{GuiGuilong}. The existence
of global weak solutions for (2CH) with $\sigma=1$ was investigated in \cite{Guan.weak}.\\
$~~~~~~$The other one is the modified two-component Camassa-Holm
system (M2CH) \cite{Holm.Naraigh}:
\begin{align}
\left\{
\begin{array}{ll}
m_t+um_x+2u_xm+\sigma\rho\overline{\rho}_x=0, \\[1ex]
\rho_t+(u\rho)_x=0,
\end{array}
\right.
\end{align}
where $m=u-u_{xx}$, $\rho=(1-\partial^2_x)(\overline{\rho}-\overline{\rho}_0)$ and $\sigma=\pm1$. Local well-posedness for (M2CH) with the initial
data in Sobolev spaces and in Besov spaces was established in \cite{Guan.Karlsen} and \cite{Kai.Yin} respectively. The blow up phenomena of strong solutions to (M2CH) were presented in \cite{Guan.Karlsen}. The existence
of global weak solutions for (M2CH) with $\sigma=1$ was investigated in \cite{Guan.weak.modified}.  The global conservative and dissipative solutions of (M2CH) were studied in \cite{Tan.Yin} and  \cite{Tan.Yin2}, respectively.\\
$~~~~~~$Recently, the authors in \cite{Wei.JDE} studied the local well-posedness and global existence of strong solutions to (1.1) under some sign condition. However, the solitons of (1.1) are not strong solutions and do not belong to the spaces $H^s(\mathbb{R})$, $s>\frac{3}{2}$. This fact motivates us to study weak solutions of (1.1). The main idea is based on the approximation of the initial data by smooth functions producing a sequence of global strong solutions $(a^n,b^n,c^n)$ of (1.1). This method was first utilized by Constantin and Molinet in \cite{Constantin.Molinet}. Due that the structure of the system (1.1) is more complex than that of the CH equation, we can not obtain the desired result under the same condition mentioned in \cite{Constantin.Molinet}. In order to obtain the existence of global weak solutions of (1.1), we have to assume that the initial data $(u_0,w_0)\in (L^1(\mathbb{R})\cap L^{1+\varepsilon}(\mathbb{R}))$, for some $\varepsilon>0$. The main difficulty is to get the uniform boundedness of $b^n$. In order to overcome this difficulty, we make good use of the special structure of the system.\\
  $~~~~~~$ The paper is organized as follows. In Section 2, we recall some properties about strong solutions of (1.1). Moveover, we give some a prior estimates which are crucial to prove our main result. In Section 3, we introduce the definition of weak solutions to (1.1) and then prove the global existence of weak solutions to (1.1).
\section{Preliminaries}
In this section we recall the global existence of strong solutions to (1.1) and some lemmas that will be used to prove our main result.
\begin{lemm}\cite{Wei.JDE}\label{Wei}
Assume that $v_0=0,~(u_0, w_0)\in (H^3(\mathbb{R}))^2$, and that $u_0=a_0-a_{0,xx}$ and $w_0=c_0-c_{0,xx}$ are nonnegative. Then the initial value problem (1.1) has a unique solution  $(u,~0,~w)\in [C(\mathbb{R}_{+};H^3(\mathbb{R}))\cap C^1(\mathbb{R}_{+};H^2(\mathbb{R}))]^3$. Moreover, $H_1(t)=\int_{\mathbb{R}}ac+a_xc_x$ and $H_2(t)=\int_{\mathbb{R}}uc_xdx=-\int_{\mathbb{R}}wa_x$ are conservation laws. For every $t\geq 0$ we have
\begin{align*}
&(1)\quad |a_x(t,x)|\leq a(t,x) \quad \text{and}\quad|c_{x}(t,x)|\leq c(t,x),\quad \forall x\in\mathbb{R}, \\
&(2)\quad u(t,x)\geq 0  \quad \text{and}\quad w(t,x)\geq 0,\quad \forall x\in\mathbb{R}, \\
&(3)\quad \|a_x(t,\cdot)\|_{L^\infty(\mathbb{R})}\leq \|a(t,\cdot)\|_{L^\infty(\mathbb{R})}\leq C\|a(t,\cdot)\|_{H^1(\mathbb{R})}\leq C\exp{[(4H_1(0)+H_2(0))t]} \quad and \\
&\|c_x(t,\cdot)\|_{L^\infty(\mathbb{R})}\leq \|c(t,\cdot)\|_{L^\infty(\mathbb{R})}
\leq C\|c(t,\cdot)\|_{H^1(\mathbb{R})}\leq C\exp{[(4H_1(0)+H_2(0))t]}, \\
&(4)\quad \|b(t,\cdot)\|_{L^\infty(\mathbb{R})}, \|b_x(t,\cdot)\|_{L^\infty(\mathbb{R})}\leq H_1(0)+\frac{1}{4}H_2(0)+\exp{[(8H_1(0)+2H_2(0))t]} .
\end{align*}
\end{lemm}
\begin{lemm}\label{Lemma2}
Assume that $v_0=0,~(u_0, w_0)\in (H^3(\mathbb{R}))^2$, and that $u_0=a_0-a_{0,xx}$ and $w_0=c_0-c_{0,xx}$ are nonnegative. And let $(u,0,w)$ be the corresponding solution to (1.1) as in Lemma \ref{Wei}. Then for any $t\in[0,T]$, there exists a constant $C$ such that
\begin{align*}
\|b(t,\cdot)\|_{H^1}\leq C(H_1(0)+H_2(0))+C\exp{[(8H_1(0)+2H_2(0))t]}.
\end{align*}
\begin{proof}
Since $v=0$, it follows from (1.1) that $4b-b_{xx}=a_{xx}c_x-c_{xx}a_x+3a_xc-3ac_x$. Note that $G_2\ast f= (4-\partial_{xx})^{-1}f$ with $G_2(x)=\frac{1}{8}e^{-2|x|}$. Applying Young's inequality, we deduce that
\begin{align*}
\|b(t,\cdot)\|_{H^1}&\leq C\int_{\mathbb{R}}|a_{xx}c_{x}-c_{xx}a_{x}+3a_xc-3ac_x|dx\\
 &\leq C\int_{\mathbb{R}}(|u(c_x+c)|+|w(a_x+a)|+|uc|+|wa|+2|a_xc|+2|ac_x|)dx\\
&\leq C\int_{\mathbb{R}}(|u(c_x+c)|+|w(a_x+a)|+|uc|+|wa|)dx+ +C\|a\|_{H^1}\|c\|_{H^1}.
\end{align*}
Thanks to Lemma \ref{Wei}, we see that $u\geq 0$, $w\geq 0$, $a_x+a\geq0$, $c_x+c\geq 0$, $a\geq0$, $c\geq0$, which leads to
\begin{align*}
\|b(t,\cdot)\|_{H^1} &\leq C\int_{\mathbb{R}} [u(c_x+c)+w(a_x+a)+uc+wa] dx+C\|a\|_{H^1}\|c\|_{H^1} \\
&\leq C(H_1(0)+H_2(0))+C\exp{[(8H_1(0)+2H_2(0))t]}.
\end{align*}
\end{proof}
\end{lemm}
Now we present some $L^p$-estimates of the strong solution to (1.1) where $ p\in[1,\infty]$.
\begin{lemm}\label{estimate}
Assume that $v_0=0,~(u_0, w_0)\in (H^3(\mathbb{R}))^2$, and that $u_0=a_0-a_{0,xx}$ and $w_0=c_0-c_{0,xx}$ are nonnegative and belong to $L^1(\mathbb{R})\cap L^{1+\varepsilon}(\mathbb{R})$ for some $\varepsilon>0$. And let $(u,0,w)$ be the corresponding solution of (1.1) as in Lemma \ref{Wei}. Then for any $t\in[0,T]$, there exists a constant $C_T$ such that
\begin{align}
\|a(t,\cdot)\|_{L^1(\mathbb{R})}=\|u(t,\cdot)\|_{L^1(\mathbb{R})}\leq e^{tC_T}\|u_0\|_{L^1(\mathbb{R})},\quad \|c(t,\cdot)\|_{L^1(\mathbb{R})}=\|w(t,\cdot)\|_{L^1(\mathbb{R})}\leq e^{tC_T}\|w_0\|_{L^1(\mathbb{R})},
\end{align}
\begin{align}
\|a(t,\cdot)\|_{L^{1+\varepsilon}(\mathbb{R})}\leq \|u(t,\cdot)\|_{L^{1+\varepsilon}(\mathbb{R})}\leq e^{tC_T}\|u_0\|_{L^{1+\varepsilon}(\mathbb{R})},\quad
\|c(t,\cdot)\|_{L^{1+\varepsilon}(\mathbb{R})}\leq \|w(t,\cdot)\|_{L^{1+\varepsilon}(\mathbb{R})}\leq e^{tC_T}\|u_0\|_{L^{1+\varepsilon}(\mathbb{R})}.
\end{align}
\begin{proof}
By density argument, we assume that $u(t,\cdot)\in C^{\infty}_0(\mathbb{R})$ and $w(t,\cdot)\in C^{\infty}_0(\mathbb{R})$. By virtue of (1.1) and integration by parts , we have
\begin{align}
&\frac{d}{dt}\int^{+\infty}_{-\infty}udx=\int^{+\infty}_{-\infty}u_tdx=\int^{+\infty}_{-\infty}u_xb+\frac{3}{2}ub_x-\frac{3}{2}u(a_xc_x-ac)dx\\
\nonumber&=\int^{+\infty}_{-\infty}\frac{1}{2}ub_x-\frac{3}{2}u(a_xc_x-ac)dx\leq \bigg\{\frac{1}{2}\|b_x\|_{L^\infty([0,T)\times\mathbb{R})}+\frac{3}{2}\|a_xc_x-ac\|_{L^\infty([0,T)\times\mathbb{R})}\bigg\}\|u\|_{L^1(\mathbb{R})}.
\end{align}
Taking advantage of Lemma \ref{Wei} and using the fact that $u\geq0$, we deduce that
\begin{align}
\|u(t,\cdot)\|_{L^1(\mathbb{R})}\leq \|u_0\|_{L^1(\mathbb{R})}+\int^t_0 C_T\|u(s,\cdot)\|_{L^1(\mathbb{R})}ds.
\end{align}
Applying Gronwall's inequality, we infer that
\begin{align}
\|u(t,\cdot)\|_{L^1(\mathbb{R})} \leq e^{tC_T}\|u_0\|_{L^1(\mathbb{R})}.
\end{align}
Since $a(t,x)$ and $u(t,x)$ are nonnegative, it follows that
\begin{align}
\|u(t,\cdot)\|_{L^1(\mathbb{R})}=\int^{+\infty}_{-\infty}u(t,x)dx=\int^{+\infty}_{-\infty}a(t,x)-a_{xx}(t,x)dx=\int^{+\infty}_{-\infty}a(t,x)dx=\|a(t,\cdot)\|_{L^1(\mathbb{R})}.
\end{align}
By the same token, we obtain
\begin{align}
\|c(t,\cdot)\|_{L^1(\mathbb{R})}=\|w(t,\cdot)\|_{L^1(\mathbb{R})}\leq e^{tC_T}\|w_0\|_{L^1(\mathbb{R})}.
\end{align}
Now we turn our attention to prove (2.2). If $\varepsilon<\infty$, By virtue of (1.1) and integration by parts , we have
\begin{align}
&\frac{d}{dt}\int^{+\infty}_{-\infty}u^{1+\varepsilon}dx=(1+\varepsilon)\int^{+\infty}_{-\infty}u^\varepsilon u_tdx=\int^{+\infty}_{-\infty}u^{1+\varepsilon}_xb+\frac{3(1+\varepsilon)}{2}u^{1+\varepsilon}b_x-\frac{3(1+\varepsilon)}{2}u^{1+\varepsilon}(a_xc_x-ac)dx\\
\nonumber&=\int^{+\infty}_{-\infty}\frac{1+3\varepsilon}{2}u^{1+\varepsilon}b_x-\frac{3(1+\varepsilon)}{2}u^{1+\varepsilon}(a_xc_x-ac)dx\\
\nonumber&\leq\bigg\{\frac{1+3\varepsilon}{2}\|b_x\|_{L^\infty([0,T)\times\mathbb{R})}+\frac{3(1+\varepsilon)}{2}\|a_xc_x-ac\|_{L^\infty([0,T)\times\mathbb{R})}\bigg\}\|u\|^{1+\varepsilon}_{L^{1+\varepsilon}(\mathbb{R})},
\end{align}
which along with $u\geq0$ leads to
\begin{align}
\frac{d}{dt}\|u\|_{L^{1+\varepsilon}(\mathbb{R})}&\leq \bigg\{\frac{1+3\varepsilon}{2(1+\varepsilon)}\|b_x\|_{L^\infty([0,T)\times\mathbb{R})}+\frac{3}{2}\|a_xc_x-ac\|_{L^\infty([0,T)\times\mathbb{R})}\bigg\}\|u\|_{L^{1+\varepsilon}(\mathbb{R})}\\
\nonumber&\leq\frac{3}{2}(\|b_x\|_{L^\infty([0,T)\times\mathbb{R})}+\|a_xc_x-ac\|_{L^\infty([0,T)\times\mathbb{R})})\|u\|_{L^{1+\varepsilon}(\mathbb{R})}.
\end{align}
Taking advantage of Lemma \ref{Wei} and Gronwall's inequality, we infer that
\begin{align}
\|u(t,\cdot)\|_{L^{1+\varepsilon}(\mathbb{R})} \leq e^{tC_T}\|u_0\|_{L^{1+\varepsilon}(\mathbb{R})}.
\end{align}
If $\varepsilon=\infty$, using a similar calculation for any $0<\delta<\infty$, and then taking limit as $\delta\rightarrow\infty$, we obtain
\begin{align}
\|u(t,\cdot)\|_{L^{\infty}(\mathbb{R})} \leq e^{tC_T}\|u_0\|_{L^{\infty}(\mathbb{R})}.
\end{align}
Note that $G_1\ast f=(1-\partial_{xx})^{-1}f$ with $G_1(x)=\frac{1}{2}e^{-|x|}$. Using Young's inequality, we deduce that
\begin{align}
\|a(t,\cdot)\|_{L^{1+\varepsilon}(\mathbb{R})}=\|G_1\ast u(t,\cdot)\|_{L^{1+\varepsilon}(\mathbb{R})}\leq \|u(t,\cdot)\|_{L^{1+\varepsilon}(\mathbb{R})}\leq e^{tC_T}\|u_0\|_{L^{\infty}(\mathbb{R})}.
\end{align}
By the same token, we get
\begin{align}
\|c(t,\cdot)\|_{L^{1+\varepsilon}(\mathbb{R})}\leq \|w(t,\cdot)\|_{L^{1+\varepsilon}(\mathbb{R})}\leq e^{tC_T}\|w_0\|_{L^{\infty}(\mathbb{R})}.
\end{align}
\end{proof}
\end{lemm}
Let us now recall a partial integration result for Bochner spaces.
\begin{lemm}\cite{Malek}\label{Bochner}
Let $T>0$. If
$$f,~g\in L^2(0,T;H^1(\mathbb{R}))\quad \text{and} \quad \frac{df}{dt}, \frac{dg}{dt}\in L^2(0,T;H^{-1}(\mathbb{R})),$$
then $f,~g$ are a.e. equal to a function continuous from $[0,T]$ into $L^2(\mathbb{R})$ and
$$\langle f(t), g(t)\rangle-\langle f(s), g(s)\rangle=\int^t_s\langle \frac{df(\tau)}{d\tau}, g(\tau)\rangle d\tau+\int^t_s\langle \frac{dg(\tau)}{d\tau}, f(\tau)\rangle d\tau$$
for all $s,t\in [0,T]$, where $\langle \cdot, \cdot\rangle$ is the $H^{-1}(\mathbb{R})$ and $H^{1}(\mathbb{R})$ duality bracket.
\end{lemm}
Throughout this paper, let $\{\rho_n\}_{n\geq1}$ denote the mollifiers
\begin{align*}
\rho_n(x)=\bigg(\int_{\mathbb{R}}\rho(y)dy\bigg)^{-1}n\rho(nx),\quad x\in\mathbb{R},\quad n\geq1,
\end{align*}
where $\rho\in C^\infty_0(\mathbb{R})$ is defined by
$$\rho(x)=\left\{\begin{array}{ll}
e^{\frac{1}{x^2-1}},~~~|x|<1, \\[1ex]
0,~~~~~~~~~|x|>1.
\end{array}
\right.$$

\section{Global weak solutions}
In this section, we first introduce the definition of weak solutions to (1.1) with $v=0$. Note that $G_1\ast f=(1-\partial_{xx})^{-1}f$ with $G_1(x)=\frac{1}{2}e^{-|x|}$. For smooth solutions of (1.1), we get
\begin{align}
a_t&=G_1\ast[u_xb+\frac{3}{2}ub_x-\frac{3}{2}u(a_xc_x-ac)]\\
\nonumber&=a_xb+\frac{1}{2}\partial_xG_1\ast(a_xb_x+a^2_xc_x+3ab-3a_xac)+\frac{1}{2}G_1\ast(ba_x+3a^2c+3aa_xc_x).
\end{align}
By the same token, we obtain
\begin{align}
c_t&=G_1\ast[w_xb+\frac{3}{2}wb_x+\frac{3}{2}w(a_xc_x-ac)]\\
\nonumber&=c_xb+\frac{1}{2}\partial_xG_1\ast(b_xc_x-a_xc^2_x+3bc+3acc_x)+\frac{1}{2}G_1\ast(bc_x-3ac^2-3a_xcc_x).
\end{align}
For simplicity, we introduce the notation
 \begin{align}
f_1=a_xb_x+a^2_xc_x+3ab-3a_xac,~~f_2=b_xc_x-a_xc^2_x+3bc+3acc_x,
\end{align}
\begin{align}
g_1=ba_x+3a^2c+3aa_xc_x,~~g_2=bc_x-3ac^2-3a_xcc_x.
\end{align}
Then (1.1) can be rewrite in the following hyperbolic type
    \begin{align}
\left\{
\begin{array}{ll}
a_{t}=a_xb+\frac{1}{2}\partial_xG_1\ast f_1+\frac{1}{2}G_1\ast g_1,\\[1ex]
c_t=c_xb+\frac{1}{2}\partial_xG_1\ast f_2+\frac{1}{2}G_1\ast g_2,\\[1ex]
4b-b_{xx}=a_{xx}c_x-c_{xx}a_x+3a_xc-3ac_x, \\[1ex]
a|_{t=0}=a_{0},~~ c|_{t=0}=c_0.\\[1ex]
\end{array}
\right.
\end{align}
\begin{defi}\label{weak}
Assume that $(a_0,c_0)\in (H^s(\mathbb{R}))^2$ with $s<\frac{5}{2}$. If $(a,c)\in [L^\infty_{loc}(0,T;H^s(\mathbb{R}))]^2$ and satisfies
\begin{align*}
\int^T_0\int_{\mathbb{R}}(a\phi_t+a_xb\phi -\frac{1}{2}G_1\ast f_1\phi_x+\frac{1}{2}G_1\ast g_1\phi) dxdt+\int_{\mathbb{R}}a_0\phi(0,x)dx=0,~~~~\forall \phi\in C^{\infty}_0((-T,T)\times\mathbb{R}),
\end{align*}
\begin{align*}
\int^T_0\int_{\mathbb{R}}(c\varphi_t+c_xb\varphi -\frac{1}{2}G_1\ast f_2\varphi_x+\frac{1}{2}G_1\ast g_2\varphi) dxdt+\int_{\mathbb{R}}c_0\varphi(0,x)dx=0,~~~~\forall \varphi\in C^{\infty}_0((-T,T)\times\mathbb{R}),
\end{align*}
\begin{align*}
\int_{\mathbb{R}}(b(t,x)\psi-b(t,x)\psi_{xx})dx = \int_{\mathbb{R}}(a_{xx}c_x-c_{xx}a_x+3a_xc-3ac_x)(t,x)\psi dx=0,~~\text{for a.e.}~~ t\in[0,T),~~\forall \psi\in C^{\infty}_0(\mathbb{R}),
\end{align*}
then $(a,b,c)$ is called a weak solution to (3.5). Moreover, if $(a(t,x),b(t,x),c(t,x))$ is a weak solution on $[0,T)$ for any $T>0$, then it is called a global weak solution to (3.5).
\end{defi}
Our main result can be stated as follow.
\begin{theo}\label{th}
Let $(a_0,c_0)\in H^1(\mathbb{R})$. Moreover $u_0=a_0-a_{0,xx}$ and $w_0=c_0-c_{0,xx}$ belong to $L^1(\mathbb{R})\cap L^{1+\varepsilon}(\mathbb{R})$ for some $\varepsilon>0$. If $u_0\geq 0$ and $w_0\geq 0$ a.e. on $\mathbb{R}$, then (3.5) has a global weak solution $(a,c)\in [W^{1,\infty}([0,T)\times\mathbb{R})\cap C([0,T);L^2({\mathbb{R}}))\cap C_w(0,T; H^1(\mathbb{R}))]^2$ for arbitrary finite $T>0$. Moreover, $(u,w)\in [L^\infty_{loc}(\mathbb{R}_+;L^1(\mathbb{R})\cap L^{1+\varepsilon}(\mathbb{R}))]^2$.
\begin{proof}
{\bf Step 1.} Without loss of generality, we assume that $\varepsilon<\infty$.  Define $a^n_0=\rho_n\ast a_0\in H^\infty(\mathbb{R})$ and $c^n_0=\rho_n\ast c_0\in H^\infty(\mathbb{R})$ for $n\geq 1$. Then we have
\begin{align}
a^n_0\rightarrow a_0  \quad\text{and}\quad  c^n_0\rightarrow c_0  \quad \text{in}\quad H^1(\mathbb{R}),\quad \text{as}\quad n\rightarrow\infty.
\end{align}
 Since $u^n_0=a^n_0-a^n_{0,xx}=\rho_n\ast u_0$ and $ w^n_0=c^n_0-c^n_{0,xx}=\rho_n\ast w_0$ for $n\geq1$, it follows that
\begin{align}
u^n_0\rightarrow u_0 \quad\text{and}\quad  w^n_0\rightarrow w_0  \quad \text{in}\quad L^1(\mathbb{R})\cap L^{1+\varepsilon}(\mathbb{R}),\quad \text{as}\quad n\rightarrow\infty.
\end{align}
Note that $u^n_0\geq 0$ and $w^n_0\geq0$. By Lemma \ref{Wei}, we obtain that there exists a global strong solution $(u^n,0,w^n)$ of (1.1) with the initial data $(u^n_0,0,w^n_0)$. Moreover $(u^n,w^n)\in [C([0,\infty);H^s(\mathbb{R}))\cap C^1([0,\infty);H^{s-1}(\mathbb{R})]^2$ for any $s\geq3$ and $u^n=a^n-a^n_{xx}\geq0$,  $w^n=c^n-c^n_{xx}\geq0$.\\
{\bf Step 2.} For fixed $T>0$, by virtue of Lemmas \ref{Wei}-\ref{Lemma2}, we have
\begin{align}
&\|a^n_x\|_{L^\infty([0,T]\times\mathbb{R})}\leq \|a^n\|_{L^\infty([0,T]\times\mathbb{R})}\leq C\|a^n\|_{L^\infty(0,T;H^1(\mathbb{R}))}\leq C\exp{[(4H^n_1(0)+H^n_2(0))T]},\\
&\|c^n_x\|_{L^\infty([0,T]\times\mathbb{R})}\leq \|c^n\|_{L^\infty([0,T]\times\mathbb{R})}\leq C\|c^n\|_{L^\infty([0,T];H^1(\mathbb{R}))}\leq C\exp{[(4H^n_1(0)+H^n_2(0))T]},\\
&\|b^n\|_{L^\infty([0,T]\times\mathbb{R})},~\|b^n_x\|_{L^\infty([0,T]\times\mathbb{R})}\leq H^n_1(0)+\frac{1}{4}H^n_2(0)+\exp{[(8H^n_1(0)+2H^n_2(0))T]}, \\
\nonumber&\|b^n\|_{L^\infty([0,T];H^1(\mathbb{R}))}\leq C\{(H^n_1(0)+H^n_2(0))+\exp{[(4H^n_1(0)+H^n_2(0))T]}\},
\end{align}
where $H^n_1(0)=\int_{\mathbb{R}}a^n_0c^n_0+a^n_{0,x}c^n_{0,x}$ and $H^n_2(0)=\int_{\mathbb{R}}u^n_0c^n_{0,x}$. Applying Cauchy-Schwarz's inequality and Young's inequality, we obtain
\begin{align}
H^n_1(0)\leq \|a^n_0\|_{H^1(\mathbb{R})}+\|c^n_0\|_{H^1(\mathbb{R})}\leq \|a_0\|_{H^1(\mathbb{R})}+\|c_0\|_{H^1(\mathbb{R})}.
\end{align}
 Since $w^n_0\geq 0$, it follows that
\begin{align}
H^n_2(0)\leq \|u^n_0\|_{L^1(\mathbb{R})}\|c^n_{0,x}\|_{L^\infty(\mathbb{R})}\leq \|u^n_0\|_{L^1(\mathbb{R})}\|c^n_0\|_{L^\infty(\mathbb{R})}\leq \|u_0\|_{L^1(\mathbb{R})}\|c_0\|_{L^\infty(\mathbb{R})}.
\end{align}
Plugging (3.11)-(3.12) into (3.8)-(3.10), we verify that $(a^n,c^n)$ is uniformly bounded in $[L^\infty(0,T;W^{1,\infty}(\mathbb{R}))\cap L^\infty(0,T;H^1(\mathbb{R}))]^2$ and $b^n$ is uniformly bounded in $L^\infty(0,T;W^{1,\infty}(\mathbb{R}))\cap L^\infty(0,T;H^1(\mathbb{R}))$.
By virtue of (3.5), we obtain
\begin{align}
a^n_t=a^n_xb^n+\frac{1}{2}\partial_xG_1\ast f^n_1+\frac{1}{2}G_1\ast g^n_1,
\end{align}
where $f^n_1=a^n_xb^n_x+(a^n_x)^2c^n_x+3a^nb^n-3a^n_xa^nc^n$ and $g^n_1=b^na^n_x+3(a^n)^2c^n+3a^na^n_xc^n_x$. \\
Since $(a^n,c^n)$ is uniformly bounded in $[L^\infty(0,T;W^{1,\infty}(\mathbb{R}))\cap L^\infty(0,T;H^1(\mathbb{R}))]^2$, it follows that $a^n_t$ is uniformly bounded in $L^\infty((0,T)\times\mathbb{R})\cap L^2((0,T)\times\mathbb{R})$. Similarly, we deduce that $c^n_t$ is uniformly bounded in $L^\infty((0,T)\times\mathbb{R})\cap L^2((0,T)\times\mathbb{R})$. Therefore, it has a subsequence such that
\begin{align}
(a^{n_k},c^{n_k})\rightharpoonup (a,c),\quad *~\text{weakly in} \quad [W^{1,\infty}((0,T)\times\mathbb{R})\cap H^1((0,T)\times\mathbb{R})]^2 \quad \text{as} \quad n_k\rightarrow\infty,
\end{align}
and
\begin{align}
(a^{n_k},c^{n_k})\xrightarrow {n_k\rightarrow\infty} (a,c), \quad \text{a.e. on} \quad (0,T)\times \mathbb{R},
\end{align}
for some $(a,c)\in [W^{1,\infty}((0,T)\times\mathbb{R})\cap H^1((0,T)\times\mathbb{R})]^2$.  By virtue of Lemma \ref{estimate}, we see that  \begin{align}
\|a^{n_k}_{xx}\|_{L^\infty(0,T;L^1(\mathbb{R})}\leq \|a^{n_k}\|_{L^\infty(0,T;L^1(\mathbb{R})}+\|u^{n_k}\|_{L^\infty(0,T;L^1(\mathbb{R})}\leq 2e^{TC_T}\|u_0\|_{L^1(\mathbb{R})}.
\end{align}
Differentiating  (3.13) with respect to $x$ yields that
\begin{align*}
a^n_{xt}=a^n_{xx}b^n+a^n_xb^n_x+\frac{1}{2}G_1\ast f^n_1-\frac{1}{2}f^n_1+\frac{1}{2}\partial_xG_1\ast g^n_1,
\end{align*}
which along with Young's inequality leads to $\|a^{n_k}_{xt}\|_{L^\infty(0,T;L^1(\mathbb{R})}\leq C_T$. Since $T<\infty$, it follows that
\begin{align}
\mathds{V}[a^{n_k}_x]=\|a^{n_k}_{xx}\|_{L^1((0,T)\times\mathbb{R})}+\|a^{n_k}_{xt}\|_{L^1((0,T)\times\mathbb{R})}\leq C_T,
\end{align}
where $\mathds{V}(f)$ is the total variation of $f\in BV([0,T]\times\mathbb{R})$. By Helly's theorem (See \cite{Natanson}), there exists a subsequence, denoted again by $a^{n_k}_x$, such that
\begin{align}
a^{n_k}_x\xrightarrow{n_k\rightarrow\infty} \alpha,  \quad \text {a.e. on} \quad (0,T)\times\mathbb{R},
\end{align}
where $\alpha\in BV((0,T)\times\mathbb{R})$ with $\mathds{V}(\alpha)\leq C_T$.
From (3.15) we have $a^{n_k}_x\xrightarrow {n_k\rightarrow\infty} a_x$ in $\mathcal{D}'((0,T)\times\mathbb{R})$. This enables us to identify $\alpha$ with $a_x$ for a.e. $t\in(0,T)\times\mathbb{R}$. Therefore
\begin{align}
a^{n_k}_x \xrightarrow{n_k\rightarrow\infty} a_x, \quad \text {a.e. on} \quad (0,T)\times\mathbb{R},
\end{align}
and $\mathds{V}(a_x)\leq C_T.$ By the same token, we deduce that
\begin{align}
c^{n_k}_x \xrightarrow{n_k\rightarrow\infty} c_x, \quad \text {a.e. on} \quad (0,T)\times\mathbb{R}.
\end{align}
~~~~Note that $G_2\ast f=(4-\partial_{xx})^{-1}f$ with $G_2(x)=\frac{1}{8}e^{-2|x|}$. By virtue of (3.5), we have
\begin{align}
b^n=G_2\ast (a^n_{xx}c^n_x-c^n_{xx}a^n_x+3a^n_xc^n-3a^nc^n_x)=G_2\ast (a^n_xw^n-c^n_xu^n+2a^n_xc^n-2a^nc^n_x).
\end{align}
By (1.1), we deduce that
\begin{align}
b^n_t&=G_2\ast(a^n_{xt}w^n-c^n_{xt}u^n)+G_2\ast(a^n_xw^n_t-c^n_xu^n_t)+2G_2\ast(a^n_{x,t}c^n-a^nc^n_{x,t})+2G_2\ast(a^n_xc^n_t-a^n_tc^n_x)\\
\nonumber&=\uppercase\expandafter{\romannumeral1}^n+\uppercase\expandafter{\romannumeral2}^n+\uppercase\expandafter{\romannumeral3}^n+\uppercase\expandafter{\romannumeral4}^n.
\end{align}
Since $a^n_x$, $a^n_t$, $c^n_x$ and $c^n_t$ are uniformly bounded in $L^\infty((0,T)\times\mathbb{R})\cap L^2((0,T)\times\mathbb{R})$, it follows from Young's inequality that {\footnote{For simplicity, we use the notation $\|\cdot\|_{L^\infty\cap L^2}$ instead of $\|\cdot\|_{L^\infty((0,T)\times\mathbb{R})\cap L^2((0,T)\times\mathbb{R})}$.}}
\begin{align}
\|\uppercase\expandafter{\romannumeral4}^n\|_{L^\infty\cap L^2}\leq 2(\|a^n_x\|_{L^\infty\cap L^2}\|c^n_t\|_{L^\infty\cap L^2}+\|c^n_x\|_{L^\infty\cap L^2}\|a^n_t\|_{L^\infty\cap L^2})\leq C_T.
\end{align}
We first consider the term $\uppercase\expandafter{\romannumeral1}^n$. By virtue of (3.5), we see that
\begin{align}
I^n&=G_2\ast[(a^n_{xx}b^n+a^n_xb^n_x+\frac{1}{2}G_1\ast f^n_1-\frac{1}{2}f^n_1+\frac{1}{2}\partial_xG_1\ast g^n_1)w^n\\
\nonumber&-(c^n_{xx}b^n+b^n_xc^n_x+\frac{1}{2}G_1\ast f^n_2-\frac{1}{2}f^n_2+\frac{1}{2}\partial_xG_2\ast g^n_2)u^n]\\
\nonumber&=G_2\ast[(a^n_{xx}b^n+\frac{1}{2}a^n_xb^n_x-\frac{1}{2}(a^n_x)^2c^n_x-\frac{3}{2}a^nb^n+\frac{3}{2}a^n_xa^nc^n+\frac{1}{2}G_1\ast f^n_1+\frac{1}{2}\partial_xG_1\ast g^n_1)w^n]\\
\nonumber&-(c^n_{xx}b^n+\frac{1}{2}c^n_xb^n_x+\frac{1}{2}a^n_x(c^n_x)^2-\frac{3}{2}b^nc^n-\frac{3}{2}a^nc^nc^n_x+\frac{1}{2}G_1\ast f^n_2+\frac{1}{2}\partial_xG_2\ast g^n_2)u^n\\
\nonumber&=G_2\ast[(a^n_{xx}c^n-c^n_{xx}a^n)b^n]+\frac{1}{2}G_2\ast[b^n_x(a^n_xw^n-c^n_xu^n)]-\frac{1}{2}G_2\ast[(a^n_xc^n_x-3a^nc^n)(a^n_xw^n+c^n_xu^n)]\\
\nonumber&+\frac{1}{2}G_2\ast[(G_1\ast f^n_1+\partial_xG_1\ast g^n_1)w^n]+\frac{1}{2}G_2\ast[(G_1\ast f^n_2+\partial_xG_1\ast g^n_2)u^n]\\
\nonumber&=\uppercase\expandafter{\romannumeral1}^n_1+\uppercase\expandafter{\romannumeral1}^n_2+\uppercase\expandafter{\romannumeral1}^n_3+\uppercase\expandafter{\romannumeral1}^n_4+\uppercase\expandafter{\romannumeral1}^n_5.
\end{align}
{\bf Bounds for $\uppercase\expandafter{\romannumeral1}^n_1$}. Since $(a^n_{xx}c^n-c^n_{xx}a^n)=(a^n_xc^n-c^n_xa^n)_x$, it follows that
\begin{align}
\uppercase\expandafter{\romannumeral1}^n_1= 4G_2\ast[(a^n_xc^n-c^n_xa^n)_xb^n]=4\partial_xG_2\ast[(a^n_xc^n-c^n_xa^n)b^n]-4G_2\ast[(a^n_xc^n-c^n_xa^n)b^n_x],
\end{align}
which together with Young's inequality leads to $\|\uppercase\expandafter{\romannumeral1}^n_1\|_{L^\infty\cap L^2}\leq C_T$.\\
{\bf Bounds for $\uppercase\expandafter{\romannumeral1}^n_2$}. In view of the fact that $a^n_xw^n-c^n_xu^n= 2a^nc^n_x-2a^n_xc^n+4b^n-b^n_{xx}$, which implies that
\begin{align}
\uppercase\expandafter{\romannumeral1}^n_2&= \frac{1}{2}G_2\ast[(2a^nc^n_x-2a^n_xc^n+4b^n-b^n_{xx})b^n_x]\\
\nonumber&=G_2\ast[(a^nc^n_x-a^n_xc^n)b^n_x]+\partial_xG_2\ast(b^n)^2-\frac{1}{4}\partial_xG_2\ast(b^n_x)^2,
\end{align}
which along with Young's inequality leads to $\|\uppercase\expandafter{\romannumeral1}^n_2\|_{L^\infty\cap L^2}\leq C_T.$\\
{\bf Bounds for $\uppercase\expandafter{\romannumeral1}^n_3$}. Thanks to $(a^n_xw^n+c^n_xu^n)=(a^nc^n-a^n_xc^n_x)_x$, we deduce that
\begin{align}
\uppercase\expandafter{\romannumeral1}^n_3&=\frac{1}{2}\partial_xG_2\ast[(a^n_xc^n_x-a^nc^n)^2]+2G_2\ast [a^nc^n(a^n_xc^n_x-a^nc^n)_x]\\
\nonumber&=\frac{1}{2}\partial_xG_2\ast[(a^n_xc^n_x-a^nc^n)^2]-2G_2\ast [(a^nc^n)_x(a^n_xc^n_x-a^nc^n)]+2\partial_xG_2\ast[a^nc^n(a^n_xc^n_x-a^nc^n)],
\end{align}
which along with Young's inequality implies that $\|\uppercase\expandafter{\romannumeral1}^n_3\|_{L^\infty\cap L^2}\leq C_T.$ \\
{\bf Bounds for $\uppercase\expandafter{\romannumeral1}^n_4$ and $\uppercase\expandafter{\romannumeral1}^n_5$}. Since $\partial_{xx}G_1\ast f= G_1\ast f-f$, it follows that
\begin{align}
\uppercase\expandafter{\romannumeral1}^n_4&=\frac{1}{2}G_2\ast[(G_1\ast f^n_1+\partial_xG_1\ast g^n_1)c^n]-\frac{1}{2}G_2\ast[(G_1\ast f^n_1+\partial_xG_1\ast g^n_1)c^n_{xx}]\\
\nonumber&=\frac{1}{2}G_2\ast[(G_1\ast f^n_1+\partial_xG_1\ast g^n_1)c^n]-\frac{1}{2}\partial_xG_2\ast[(G_1\ast f^n_1+\partial_xG_1\ast g^n_1)c^n_{x}]\\
\nonumber&+\frac{1}{2}G_2\ast[(\partial_xG_1\ast f^n_1+G_1\ast g^n_1-g^n_1)c^n_{x}].
\end{align}
Note that $f^n_1$ and $g^n_1$ are bounded in $L^\infty((0,T)\times\mathbb{R})\cap L^2((0,T)\times\mathbb{R})$. Taking advantage of Young's inequality yields that  $\|\uppercase\expandafter{\romannumeral1}^n_4\|_{L^\infty\cap L^2}\leq C_T.$ By the same token, we have $\|\uppercase\expandafter{\romannumeral1}^n_5\|_{L^\infty\cap L^2}\leq C_T.$ Thus, we show that $\|\uppercase\expandafter{\romannumeral1}^n\|_{L^\infty\cap L^2}\leq C_T.$ Since the estimate for $\uppercase\expandafter{\romannumeral3}^n$ is similar to that of $\uppercase\expandafter{\romannumeral1}^n$, we omit the details here.\\
Now we turn our attention to estimate the term $\uppercase\expandafter{\romannumeral2}^n$. By virtue of (3.5), we have
\begin{align}
\uppercase\expandafter{\romannumeral2}^n&=G_2\ast\{a^n_x[w^n_xb^n+\frac{3}{2}w^nb^n_x+\frac{3}{2}w^n(a^n_xc^n_x-a^nc^n)]-c^n_x[u^n_xb^n+\frac{3}{2}u^nb^n_x-\frac{3}{2}u^n(a^n_xc^n_x-a^nc^n)]\}\\
\nonumber&=G_2\ast[(a^n_xw^n_x-c^n_xu^n_x)b^n]+\frac{3}{2}G_2\ast[(a^n_xw^n-c^n_xu^n)b^n_x]+\frac{3}{2}G_2\ast[(w^na^n_x+c^n_xu^n)(a^n_xc^n_x-a^nc^n)]\\
\nonumber&=\uppercase\expandafter{\romannumeral2}^n_1+\uppercase\expandafter{\romannumeral2}^n_2+\uppercase\expandafter{\romannumeral2}^n_3.
\end{align}
Since $a^n_xw^n_x-c^n_xu^n_x=a^n_{xxx}c^n_x-a^n_xc^n_{xxx}=(a^n_{xx}c^n_x-a^n_xc^n_{xx})_x=(3a^nc^n_x-3a^n_xc^n+4b^n-b^n_{xx})_x$, it follows that
\begin{align}
\uppercase\expandafter{\romannumeral2}^n_1&=\frac{1}{2}G_2\ast[(3a^nc^n_x-3a^n_xc^n+4b^n-b^n_{xx})_xb^n]\\
\nonumber&=\frac{3}{2}\partial_xG_2\ast[(a^nc^n_x-a^n_xc^n)b^n]-\frac{3}{2}G_2\ast[(a^nc^n_x-a^n_xc^n)b^n_x]+\partial_xG_2\ast(b^n)^2-\frac{1}{2}G_2\ast (b^n_{xxx}b^n).
\end{align}
From the above identity, it is sufficient to bound for $G_2\ast (b^n_{xxx}b^n)$. Indeed,
\begin{align}
G_2\ast (b^n_{xxx}b^n)&=\partial_xG_2\ast(b^n_{xx}b^n)-G_2\ast(b^n_{xx}b^n_x)\\
\nonumber&=\partial_{xx}G_2\ast(b^n_{x}b^n)-\partial_xG_2\ast(b^n_x)^2-\frac{1}{2}\partial_xG_2\ast(b^n_x)^2\\
\nonumber&=4G_2\ast(b^n_{x}b^n)-b^n_xb^n-\frac{3}{2}\partial_xG_2\ast(b^n_x)^2.
\end{align}
By virtue of Young's inequality, we obtain $\|\uppercase\expandafter{\romannumeral2}^n_1\|_{L^\infty\cap L^2}\leq C_T.$
By the similar estimates as for $\uppercase\expandafter{\romannumeral1}^n_2$ and $\uppercase\expandafter{\romannumeral1}^n_3$, we infer that  $\|\uppercase\expandafter{\romannumeral2}^n_2\|_{L^\infty\cap L^2},~\|\uppercase\expandafter{\romannumeral2}^n_3\|_{L^\infty\cap L^2}\leq C_T$.
   From the above argument, we prove that $b^n_t$ is bounded in $L^\infty((0,T)\times\mathbb{R})\cap L^2((0,T)\times\mathbb{R}) $. Moreover, there exists a subsequence such that
\begin{align}
b^{n_k}\rightharpoonup b,\quad *~\text{weakly in} \quad W^{1,\infty}((0,T)\times\mathbb{R})\cap H^1((0,T)\times\mathbb{R}) \quad \text{as} \quad n_k\rightarrow\infty,
\end{align}
and
\begin{align}
b^{n_k}\xrightarrow {n_k\rightarrow\infty} b, \quad \text{a.e. on} \quad (0,T)\times \mathbb{R},
\end{align}
for some $b\in W^{1,\infty}((0,T)\times\mathbb{R})\cap H^1((0,T)\times\mathbb{R})$. \\
By virtue of Young's inequality, we have
\begin{align}
\|b^{n_k}_{xx}\|_{L^1((0,T)\times\mathbb{R})}\leq \|b^n(t,\cdot)\|_{L^1((0,T)\times\mathbb{R})}+\|(a^n_{xx}c^n_x-c^n_{xx}a^n_x+3a^n_xc^n-3a^nc^n_x)\|_{L^1((0,T)\times\mathbb{R})}\leq C_T.
\end{align}
By differentiating both sides of (3.22) with respect to $x$, we obtain that
\begin{align}
b^{n_k}_{t,x}=\uppercase\expandafter{\romannumeral1}^{n_k}_x+\uppercase\expandafter{\romannumeral2}^{n_k}_x+\uppercase\expandafter{\romannumeral3}^{n_k}_x+\uppercase\expandafter{\romannumeral4}^{n_k}_x.
\end{align}
Thanks to $\partial_xG_2\in L^p$ for any $1\leq p\leq \infty$, one can follows the similar proof as bound for $b^n_t$ to deduce that
\begin{align}
\|b^{n_k}_{tx}\|_{L^1((0,T)\times\mathbb{R})}\leq C_T.
\end{align}
By the same token as $a^{n_k}_x$, we deduce that there exists a subsequence denoted again by $b^{n_k}$, such that
  \begin{align}
b^{n_k}_x \xrightarrow{n_k\rightarrow\infty} b_x, \quad \text {a.e. on} \quad (0,T)\times\mathbb{R},
\end{align}
and $\mathds{V}(b_x)\leq C_T.$
For any fixed $t\in(0,T)$, we have $f^n_1,~f^n_2,~g^n_1,~g^n_2$ are uniformly bounded in $L^\infty(\mathbb{R})$. Therefore,  there exists a subsequence such that
\begin{align}
(f^{n_k}_1(t,\cdot),f^{n_k}_2(t,\cdot), g^{n_k}_1(t,\cdot), g^{n_k}_2(t,\cdot))\rightharpoonup (\widetilde{f_1}(t,\cdot),\widetilde{f_2}(t,\cdot),\widetilde{g_1}(t,\cdot),\widetilde{g_2}(t,\cdot)),\quad *~\text{weakly in} ~~ [L^\infty(\mathbb{R})]^4 ~~ \text{as}~~ n_k\rightarrow\infty.
\end{align}
By virtue of (3.19), (3.20) and (3.37), we deduce that $(\widetilde{f_1},\widetilde{f_2},\widetilde{g_1},\widetilde{g_2})=(f_1,f_2,g_1,g_2)$ for a.e. $t\in(0,T)$. Since $G_1(x)\in L^1(\mathbb{R})$, it follows that
\begin{align}
(G_1\ast f^{n_k}_1, G_1\ast f^{n_k}_2, G_1\ast g^{n_k}_1, G_1\ast g^{n_k}_2)\xrightarrow{n_k\rightarrow\infty}(G_1\ast f_1, G_1\ast f_2, G_1\ast g_1, G_1\ast g_2).
\end{align}
Noticing that $(a^n,b^n,c^n)\in [C^1((0,T);H^\infty)]^3$ is the strong solution of (3.5), we have
\begin{align}
\int^T_0\int_{\mathbb{R}}(a^{n_k}\phi_t+a^{n_k}_xb^{n_k}\phi -\frac{1}{2}G_1\ast f^{n_k}_1\phi_x+\frac{1}{2}G_1\ast g^{n_k}_1\phi) dxdt+\int_{\mathbb{R}}a^{n_k}_0\phi(0,x)dx=0,~~~~\forall \phi\in C^{\infty}_0((-T,T)\times\mathbb{R}),
\end{align}
\begin{align}
\int^T_0\int_{\mathbb{R}}(c^{n_k}\varphi_t+c^{n_k}_xb^{n_k}\varphi -\frac{1}{2}G_1\ast f^{n_k}_2\varphi_x+\frac{1}{2}G_1\ast g^{n_k}_2\varphi) dxdt+\int_{\mathbb{R}}c^{n_k}_0\varphi(0,x)dx=0,~~~~\forall \varphi\in C^{\infty}_0((-T,T)\times\mathbb{R}),
\end{align}
Taking limit as $n_k\rightarrow \infty$ in the above identities, we obtain
\begin{align*}
\int^T_0\int_{\mathbb{R}}(a\phi_t+a_xb\phi -\frac{1}{2}G_1\ast f_1\phi_x+\frac{1}{2}G_1\ast g_1\phi) dxdt+\int_{\mathbb{R}}a_0\phi(0,x)dx=0,~~~~\forall \phi\in C^{\infty}_0((-T,T)\times\mathbb{R}),
\end{align*}
\begin{align*}
\int^T_0\int_{\mathbb{R}}(c\varphi_t+c_xb\varphi -\frac{1}{2}G_1\ast f_2\varphi_x+\frac{1}{2}G_1\ast g_2\varphi) dxdt+\int_{\mathbb{R}}c_0\varphi(0,x)dx=0,~~~~\forall \varphi\in C^{\infty}_0((-T,T)\times\mathbb{R}).
\end{align*}
{\bf Step 3.} According to Definition \ref{weak}, it is sufficient to prove that
$(a,b,c)$ satisfies that
\begin{align*}
\int_{\mathbb{R}}(b(t,x)\psi-b(t,x)\psi_{xx})dx = \int_{\mathbb{R}}(a_{xx}c_x-c_{xx}a_x+3a_xc-3ac_x)(t,x)\psi dx=0,~~\text{for a.e.}~~ t\in[0,T),~~\forall \psi\in C^{\infty}_0(\mathbb{R}).
\end{align*}
By virtue of  (3.15), (3.19), (3.20) and (3.33), we deduce that
\begin{align}
\int_{\mathbb{R}}(b^{n_k}(t,x)\psi-b^{n_k}(t,x)\psi_{xx})dx\xrightarrow{n_k\rightarrow \infty} \int_{\mathbb{R}}(b(t,x)\psi-b(t,x)\psi_{xx})dx ~~\text{for a.e.}~~ t\in[0,T),~~\forall \psi\in C^{\infty}_0(\mathbb{R}),
\end{align}
\begin{align}
\int_{\mathbb{R}}(3a^{n_k}_xc^{n_k}-3a^{n_k}c^{n_k}_x)(t,x)\psi dx\xrightarrow{n_k\rightarrow \infty} \int_{\mathbb{R}}(3a_xc-3ac_x)(t,x)\psi dx ~~\text{for a.e.}~~ t\in[0,T),~~\forall \psi\in C^{\infty}_0(\mathbb{R}).
\end{align}
For any fixed $t\in(0,T)$, taking advantage of Lemma \ref{estimate}, we have
\begin{align}
\|u^{n_k}(t,\cdot)\|_{L^{1+\varepsilon}}\leq C_T\|u^{n_k}_0\|_{L^{1+\varepsilon}}\leq C_T\|u_0\|_{L^{1+\varepsilon}}, \quad \|w^{n_k}(t,\cdot)\|_{L^{1+\varepsilon}}\leq C_T\|w^{n_k}_0\|_{L^{1+\varepsilon}}\leq C_T\|w_0\|_{L^{1+\varepsilon}},
\end{align}
which along with Young's inequality lead to
\begin{align}
\|a^{n_k}_{xx}(t,\cdot)\|_{L^{1+\varepsilon}}\leq C_T, \quad \|c^{n_k}_{xx}(t,\cdot)\|_{L^{1+\varepsilon}}\leq C_T.
\end{align}
Therefore there exists a subsequence, denoted again by $(a^{n_k}_{xx}(t,\cdot),c^{n_k}_{xx}(t,\cdot))$, such that
\begin{align}
(a^{n_k}_{xx}(t,\cdot),c^{n_k}_{xx}(t,\cdot))\rightharpoonup (a_{xx}(t,\cdot),c_{xx}(t,\cdot)) \quad \text{in} \quad L^{1+\varepsilon}(\mathbb{R}).
\end{align}
Since $W^{2,1+\varepsilon}_{loc}(\mathbb{R})\hookrightarrow\hookrightarrow W^{1,\infty}_{loc}(\mathbb{R})$, it follows that
\begin{align}
(a^{n_k}_x(t,\cdot),c^{n_k}_x(t,\cdot))\xrightarrow{n_k\rightarrow\infty} (a_x(t,\cdot),c_x(t,\cdot)) \quad \text {in} \quad  L^\infty_{loc}(\mathbb{R}).
\end{align}
For any  $\psi\in C^{\infty}_0(\mathbb{R})$, we have
\begin{align}
\int_{\mathbb{R}}(a^{n_k}_{xx}c^{n_k}_x-a_{xx}c_x)\psi dx
=\int_{\mathbb{R}}(a^{n_k}_{xx}-a_{xx})c_{x}\psi dx+\int_{\mathbb{R}}a^{n_k}_{xx}(c^{n_k}_x-c_x)\psi dx.
\end{align}
Using the fact that $\|c_x(t,\cdot)\|_{L^\infty(\mathbb{R})}\leq \liminf_{n_k\rightarrow\infty}\|c^{n_k}_x(t,\cdot)\|_{L^\infty(\mathbb{R})}\leq C_T$ and by virtue of (3.48), we deduce that
\begin{align}
\lim_{n_k\rightarrow\infty}\int_{\mathbb{R}}(a^{n_k}_{xx}-a_{xx})c_{x}\psi dx=0.
\end{align}
Suppose that $Supp~{\psi}\subseteq (-K,K)$ with $K\geq 0$. Then, we see that
\begin{align}
\int_{\mathbb{R}}a^{n_k}_{xx}(c^{n_k}_x-c_x)\psi dx&= \int^K_{-K}a^{n_k}_{xx}(c^{n_k}_x-c_x)\psi dx\leq \|a^{n_k}_{xx}(t,\cdot)\|_{L^1(\mathbb{R})}\|c^{n_k}_x-c_x\|_{L^\infty(-K,K)}\|\psi\|_{L^\infty}\\
\nonumber&\leq C_T\|c^{n_k}_x-c_x\|_{L^\infty(-K,K)}\rightarrow 0  \quad \text{as} \quad n_k\rightarrow\infty.
\end{align}
Taking the limit as $n_k\rightarrow \infty$ in (3.50), we get
$\lim_{n_k\rightarrow\infty}\int_{\mathbb{R}}(a^{n_k}_{xx}c^{n_k}_x-a_{xx}c_x)\psi dx=0.$ By the same token we have that $\lim_{n_k\rightarrow\infty}\int_{\mathbb{R}}(c^{n_k}_{xx}a^{n_k}_x-c_{xx}a_x)\psi dx=0.$  Since that $T$ can be taken arbitrarily, we show that $(a,b,c)$ is indeed a global weak solution of (3.5) and belongs to $[W^{1,\infty}((0,T)\times\mathbb{R})]^3$. \\
{\bf Step 4.} Note that $(\partial_t a^{n_k}(t,\cdot),\partial_t b^{n_k}(t,\cdot), \partial_t c^{n_k}(t,\cdot))$ is uniformly bounded in $L^2(\mathbb{R})$ for any $t\in(0,T)$. Hence, the map $t\mapsto (a^{n_k}(t,\cdot),b^{n_k}(t,\cdot), c^{n_k}(t,\cdot))\in (H^1(\mathbb{R})^3)$ is weakly equicontinuous on $[0,T]$. It follows from the Arzela-Ascoli theorem that $(a^{n_k}(t,\cdot),b^{n_k}(t,\cdot), c^{n_k}(t,\cdot))$ contains a subsequence, denoted again by  $(a^{n_k}(t,\cdot),b^{n_k}(t,\cdot), c^{n_k}(t,\cdot))$ converges weakly in $[H^1(\mathbb{R})]^3$ uniformly in $t$.  The limit function $(a,b,c)\in [C_w([0,T);H^1(\mathbb{R}))]^3$. \\
By virtue of Fatou's Lemma, we have
\begin{align}
\|a_t(t,\cdot)\|_{L^\infty_T(L^2(\mathbb{R}))}\leq \liminf_{n_k\rightarrow\infty}\|a^{n_k}(t,\cdot)\|_{L^\infty_T(L^2(\mathbb{R}))} \leq C_T,
\end{align}
\begin{align}
\|c_t(t,\cdot)\|_{L^\infty_T(L^2(\mathbb{R}))}\leq \liminf_{n_k\rightarrow\infty}\|c^{n_k}(t,\cdot)\|_{L^\infty_T(L^2(\mathbb{R}))} \leq C_T.
\end{align}
Taking advantage of Lemma \ref{Bochner}, we see that $(a,c)\in C([0,T);L^2(\mathbb{R}))$.
\end{proof}
\end{theo}
\begin{rema} By virtue of Lemma \ref{Wei}, we have used the conservation law $H_2(t)=\int_{\mathbb{R}}uc_xdx=\int_{\mathbb{R}}u_0c'_0dx$ to obtain the desired estimates, which implies that $u_0$ at least belongs to $L^1(\mathbb{R})$. However, the additional condition $(u_0, w_0)\in L^{1+\varepsilon}(\mathbb{R})$ in Theorem \ref{th} is technical and unnatural. How to get rid of this condition is still an open problem.
\end{rema}
\begin{rema}
In view of an interpolation argument, one can obtain that the solution $(a,c)$ of (3.5) belongs to $C([0,T)\times \mathbb{R})$ for arbitrary finite $T>0$.
\end{rema}
\begin{rema}
The condition $u_0\geq0$ and $w_0\geq 0$ in Theorem \ref{th} can be replaced by  $u_0$ and $w_0$ don't change sign. One can follow the similar step to get the global existence of weak solution to (3.5).
\end{rema}

{\bf Acknowledgements}.  This work was
partially supported by NNSFC (No.11271382), RFDP (No.
20120171110014), the Macao Science and Technology Development Fund (No. 098/2013/A3) and the key project of Sun Yat-sen University.

\end{document}